\documentclass[12pt]{article}
\usepackage{amssymb}
\usepackage{amsfonts}
\usepackage{amsmath}
\usepackage[usenames]{color}
\usepackage{mathrsfs}
\usepackage{amsfonts}
\usepackage{amssymb,amsmath}
\usepackage{CJK}
\usepackage{cite}
\usepackage{cases}
\usepackage{amsthm}\usepackage{multirow}
\usepackage{graphicx}
\usepackage{float}
\usepackage{makecell}
\usepackage{latexsym}
\usepackage{CJK}

\def\cl{\centerline}

\def\vs{\vspace*}

\def\ni{\noindent}

\numberwithin{equation}{section}
\newtheorem{theo}{Theorem}[section]
\newtheorem{defi}[theo]{Definition}
\newtheorem{coro}[theo]{Corollary}
\newtheorem{lemm}[theo]{Lemma}

\newtheorem{prop}[theo]{Proposition}

\newtheorem{rema}[theo]{Remark}

\begin{document}
\begin{center}
\cl{\large\bf \vs{8pt} Hypersurfaces with Constant Mean
Curvature on Finsler manifolds\,{$^*\,$}}
\footnote {$^*\,$ Project supported by NNSFC(No. 12261034).

}
\cl{Yali Chen$^1$, Qun He$^2$$^\dag\,$ and Yantong Qian$^{2}$}

\cl{\small 1 School of Mathematics and Statistics, Anhui Normal University, Wuhu, Anhui, 241000, China.}
\cl{\small
2 School of Mathematical Sciences, Tongji University, Shanghai,
200092, China.}
\end{center}

{\small
\parskip .005 truein
\baselineskip 3pt \lineskip 3pt

\noindent{{\bf Abstract:}
In this paper, we give the geometric meaning of hypersurfaces with constant mean curvature in a Finsler manifold by using volume preserving variation. Then we give the correspondence between the mean curvature of submanifolds by homothetic navigation, which means that some geometric properties of submanifolds are the same. Finally, we deduce a Heintze-Karcher type inequality and prove an Alexandrov type theorem in special Finsler spaces.
\vs{5pt}

\ni{\bf Key words:}
hypersurface; mean curvature; volume preservation variation; homothetic navigation.}

\ni{\it Mathematics Subject Classification (2010):} 53C20, 53C60, 53C42.}
\parskip .001 truein\baselineskip 6pt \lineskip 6pt
\section{Introduction}
In Riemannian geometry, many geometers are interested in hypersurfaces with constant anisotropic mean curvatures. Andrews studied convex hypersurfaces in Euclidean spaces evolving by anisotropic analogues of the volume-preserving mean curvature flow \cite{BA}. Koiso and Palmer obtained the stability of surfaces with constant anisotropic mean curvature \cite{MP}. Ma and Xiong applied the evolution method to present a new proof of the Alexandrov type theorem for constant anisotropic mean curvature hypersurfaces in Euclidean spaces \cite{MX}. Jorge et al. formulated a variational notion of anisotropic mean curvature for immersed hypersurfaces in arbitrary Riemannian manifolds and proved that hypersurfaces with constant anisotropic mean curvature are characterized as critical points of an elliptic parametric functional subject to a volume constraint \cite{JM}.

Compared with Riemannian geometry, the study on constant mean curvature is more difficult in Finsler geometry. There are many different ways to define mean curvature on Finsler submanifolds \cite{SS,S,HS,BYW}. Later, He et al. introduced the notion of isoparametric hypersurface in Finsler spaces and gave the relation between hypersurfaces with constant mean curvature and isoparametric hypersurfaces \cite{HYS}.

Navigation is an important technique to study geometric properties in Finsler geometry. He and Dong studied isoparametric hypersurfaces in a Randers space form and a Funk-type space using navigation technique and got complete classifications of isoparametric hypersurfaces in Randers space forms \cite{HD,HDY}. Xu et al. discussed the correspondences of Jacobi fields and isoparametric hypersurfaces by homothetic navigation \cite{XM1}.

In this paper, we mainly study the mean curvature of a Finsler submanifold. We give a unified definition of mean curvature and volume preserving variation with respect to any induced volume form. Furthermore, we study the relationship between the principal curvatures by homothetic navigation and deduce the Heintz-Karcher type inequality. The main results are given as follows,

\begin{theo} \label{thm1_1}
Let $\phi : M\rightarrow N$ be an embedded hypersurface in a Finsler manifold $(N,F)$. Then $\phi$ is a critical point of any volume preserving variation if and only if the mean curvature of $M$ is a non-zero constant in $N$ with respect to any induced volume.
\end{theo}
\begin{theo} \label{thm1_2}
Let $M$ be an oriented hypersurface in an $n$-dimensional Finsler manifold $(N,F)$ and $\tilde{F}$ be the Finsler metric on $N$ determined by the navigation datum $(F,W)$. If $W$ is a homothetic vector field, then the anisotropic mean curvature of $M$ with respect to $F$ is constant if and only if it is constant with respect to $\tilde{F}$. Especially when $W$ is a Killing vector field, the $r$-th mean curvature of $M$ with respect to $F$ is constant for some $r = 1,\ldots, n-1$ if and only if that with respect to $\tilde{F}$ is constant.
\end{theo}
\begin{theo}\label{thm1_3}
Let $(N,\tilde{F})$ be an $n$-dimensional Finsler space determined by the navigation datum $(F,W)$, where $F$ is the Minkowski metric and $W$ is a homothetic vector field with dilation $c$. Let $M=\partial U$ be a closed oriented hypersurface in $(N,\tilde{F})$, where $U$ is a compact domain. If the anisotropic mean curvature $\hat{\tilde{H}}_{\tilde{\xi}}$ of $M$ with respect to the inner normal vector $\tilde{\xi}$ in $(N,\tilde{F})$ satisfies $\hat{\tilde{H}}_{\tilde{\xi}}>c$, then we have
\begin{align}\label{1111}
\int_{M}\frac{1}{\hat{\tilde{H}}_{\tilde{\xi}}-c}d\mu_{\tilde{\xi}}\geq nV(U),
\end{align}
where $d\mu_{\tilde{\xi}}$ is the induced volume form by $\tilde{\xi}$ from Busemann-Hausdorff volume of $(N^{n},\tilde{F})$ and $V(U)$ is the volume of $U$. Moreover, the equality holds if and only if $M$ is anisotropic umbilical with respect to $F$ and $\tilde{F}$, in this case, $M$ must be a Minkowski hypersphere.
\end{theo}
\begin{theo}\label{thm1_4}
Let $(N,\tilde{F})$ be an $n$-dimensional Finsler space determined by the navigation datum $(F,W)$ and $M$ be a closed oriented hypersurface in $N$, where $F$ is a Minkowski metric and $W$ is a homothetic vector field. If the anisotropic mean curvature of $M$ is constant in $(N,\tilde{F})$ , then $M$ must be a Minkowski hypersphere.  Especially when $W$ is a Killing vector field, if the $r$-th mean curvature of $M$ with respect to $\tilde{F}$ is constant for some $r = 1,\ldots, n-1$, then $M$ must be a Minkowski hypersphere.
\end{theo}
\begin{rema}\label{rema1}
From \cite{HMFL}, $(N,\tilde{F})$ in both Theorem \ref{thm1_3} and \ref{thm1_4} have non-positive constant flag curvature. The anisotropic mean curvature and $r$-th mean curvature in Theorem \ref{thm1_2}$\sim$\ref{thm1_4} are defined by principal curvatures with respect to a given normal vector field $\xi$. For the relation of the anisotropic mean curvature and the mean curvature see Lemma \ref{lem0_11}. For the hypersurface $M$ in $(N,F)$, there exist two global normal vector fields $\xi_{\pm}$, which point to the two sides of $M$. In general, $\xi_{-}\neq-\xi_{+}$ and the Minkowski hyperspheres are different if the orientations of $M$ are different.
\end{rema}

The contents of this paper is organized as follows. In section 2, some fundamental concepts and formulas which are necessary for the present paper are given. In section 3, we introduce the volume preserving variation and give the proof of Theorem \ref{thm1_1}. In section 4, we give the relationship of principal curvature by homothetic navigation and prove Theorem \ref{thm1_2}. In section 5, we
give a Heintz-Karcher type inequality and an Alexandrov type theorem in a
Minkowski space and give the proofs of Theorem \ref{thm1_3} and Theorem \ref{thm1_4}.
\section{Preliminaries}
\subsection{Finsler manifolds}
Let~$(N,F)$ be an~$n$-dimensional oriented smooth Finsler manifold and~$TN$ be the tangent bundle over~$N$ with local coordinates
$(x,y)$, where~$x=(x^1,\cdots ,x^n)$ and~$y=(y^1,\cdots ,y^n)$. The fundamental form~$g$ of~$(N,F)$ is given by
$$g=g_{ij}(x,y)dx^{i} \otimes dx^{j}, ~~~~~~~g_{ij}(x,y)=\frac{1}{2}[F^{2}] _{y^{i}y^{j}}.$$

The projection~$\pi: TN\rightarrow N$ gives rise to the pull-back bundle~$\pi^{\ast}TN$ and its dual bundle
$\pi^{\ast}T^{\ast}N$ over~$TN\backslash\{0\}$. We shall work on $TN\backslash\{0\}$ and rigidly use only objects that are invariant under positive rescaling in $y$, so that one may view them as objects on the projective sphere bundle $SN$ using homogeneous coordinates. There exists the unique
\emph{Chern connection}~$\nabla$ with~$\nabla
\frac{\partial}{\partial x^i}=\omega_{i}^{j}\frac{\partial}{\partial
x^{j}}=\Gamma^{i}_{jk}dx^k\otimes\frac{\partial}{\partial x^{j}}$ on the pull-back bundle~$\pi^{\ast}TN$.

The \emph{$\mathbf{S}$-curvature} of $(N,F)$ with respect to the volume form $d\mu=\sigma(x)dx$ is defined by \begin{align}\label{1.4.19}
\mathbf{S}(x,y)=\frac{\partial G^{i}}{\partial y^{i}}-y^{i}\frac{\partial}{\partial x^{i}}(\ln \sigma(x)),
\end{align}
where $\sigma(x)$ is a volume measure function, $dx=dx^{1}\wedge\cdots\wedge dx^{n}$ and
\begin{align*}
G^{i}=\frac{1}{4}g^{il}\left\{[F^{2}]_{x^{k}y^{l}}y^{k}-[F^{2}]_{x^{l}}\right\}
\end{align*}
are the \emph{geodesic coefficients}.

For~$X=X^{i}\frac{\partial}{\partial x^{i}}\in \Gamma(TN)$, the \emph{covariant derivative} of~$X$ along
$v=v^i\frac{\partial}{\partial x^{i}}\in T_{x}N$ with respect to a reference vector~$w\in T_{x}N\backslash 0$ is defined by
\begin{align}
{\nabla}^{w}_{v}X(x):&=\left\{v^{j}\frac{\partial X^{i}}{\partial x^{j}}(x)+{ \Gamma}^{i}_{jk}(w)v^{j}X^{k}(x)\right\}\frac
{\partial}{\partial x^{i}},\label{Z1}
\end{align}
 where
$\Gamma^{i}_{jk}$ is the connection coefficients of \emph{Chern connection}.

If $\gamma=\gamma(t)$ satisfies $\nabla^{\dot{\gamma}}_{\dot{\gamma}}\dot{\gamma}\equiv0$, then $\gamma=\gamma(t)$ is called a \textit{geodesic} of $F$. If a vector field $J(t)$ along $\gamma(t)$ satisfies
$$\nabla_{\dot{\gamma}}^{\dot{\gamma}}\nabla_{\dot{\gamma}}^{\dot{\gamma}}J+\textbf{R}_{\dot{\gamma}}(J)=0,$$
then $J(t)$ is called a \textit{Jacobi field} along $\gamma(t)$, where $\textbf{R}_{y}$ is the \textit{flag curvature tensor} and
$$R^{i}_{k}(y)=2\frac{\partial G^{i}}{\partial x^{k}}-y^{j}\frac{\partial G^{i}}{\partial x^{j}\partial y^{k}}+2G^{j}\frac{\partial G^{i}}{\partial y^{j}\partial y^{k}}-\frac{\partial G^{i}}{\partial y^{j}}\frac{\partial G^{j}}{\partial y^{k}}.$$
Set $P=\textmd{span}\{y,u\}$, the flag curvature $K(P;y)$ of $(N,F)$ can be defined by $K(P;y)=F^{-2}g_{y}(\textbf{R}_{y}(u),u)$.

Let~${\mathcal L}:TN\rightarrow T^{\ast}N$ denote the \emph{Legendre transformation}, which satisfies~${\mathcal L}(\lambda
y)=\lambda {\mathcal L}(y)$ for all~$\lambda>0,~y\in T_{x}N$ and
\begin{align}\mathcal L(y)&=F(y)[F]_{y^{i}}(y)dx^i,~~\forall y\in T_{x}N\setminus \{0\},\\
\mathcal L^{-1}(\tau)&=F^*(\tau)[F^*]_{\tau_{i}}(\tau)\frac{\partial}{\partial x^i},~~\forall \tau\in T_{x}^*N\setminus \{0\},\label{Z01}\end{align}
 where~$F^*$ is the dual metric of~$F$.
 In general, $F^*(-\tau)\neq F^*(\tau)$ and $\mathcal L^{-1}(-\tau)\neq-\mathcal L^{-1}(\tau)$. So for any $\xi\in S_xN=\{X\in T_xN~|~F(X)=1\}$, we denote $$\nu=\mathcal L\xi,~~\xi_{-}=\frac{{\mathcal L}^{-1}(-\nu)}{F^{*}(-\nu)}.$$
For a smooth function~$f: N\rightarrow \mathbb{R}$, the \emph{gradient vector} of~$f$ at~$x$ is defined as~$\nabla f(x)={\mathcal
L}^{-1}(df(x))\in T_{x}N$.
\subsection{Immersed Submanifolds}
Let $(N,F)$ be an $n$-dimensional Finsler manifold and $\phi: M\to (N,F)$ be an $m$-dimensional immersion. Here and from now on, we will use the following convention of index ranges unless otherwise stated,
$$1\leq i, j, \cdots\leq n,\ \ \ \ 1\leq a,b,\cdots\leq m<n.$$
For simplicity, we will denote $\phi(x)$ and $d\phi(X)$ by $x$ and $X$, respectively. Let
\begin{equation}\label{ss1}
\mathcal{V}(M)=\{(x,\tau)~|~x\in M,\tau\in T_{x}^{*}N,\tau (X)=0,\forall X\in T_xM\},
\end{equation}
be the {\it normal bundle} of $\phi$ or $M$ and $\mathcal{N}M={\mathcal L}^{-1}(\mathcal{V}(M))\subset TN$ \cite{SS}. Moreover, we denote the {\it unit normal bundle} of $M$ by
$$\mathcal{V}^0(M)=\{\nu\in \mathcal{V}(M)|~F^*(\nu)=1 \},$$ and $ \mathcal{N}^{0}M={\mathcal L}^{-1}(\mathcal{V}^{0}(M))=\{\xi\in\mathcal{N}(M)~|~F(\xi)=1\}.$ We call $\xi\in \mathcal{\Gamma}(\mathcal{N}^{0}M)$ the \emph{unit normal vector field} of $M$.

The metrics on submanifolds can be induced in two ways in Finsler geometry. One is to induce the metric by an isometric immersion. In this case, $\overline{F}=\phi^{*}F$ is also a Finsler metric. Furthermore, we have
$$\bar{g}_{ab}(x,y)=\frac{1}{2}[\bar{F}^{2}]_{y^{a}y^{b}}=g_{ij}(\tilde x,\tilde y)\phi^{i}_{a}\phi^{j}_{b},$$
where
$$\tilde x^{i}=\phi^{i}(x), \tilde y^{i}=\phi^{i}_{a}y^{a}, \phi^{i}_{a}=\frac{\partial \phi^{i}}{\partial x^{a}}.$$
Set
\begin{align}\label{111}
h^{i}=\phi^{i}_{ab}y^{a}y^{b}-2\phi^{i}_{a}\bar{G}^{a}+2G^{i},
\end{align}
$$h(y) := \frac{h^{i}(y)}{\bar{F}^{2}(y)}\frac{\partial}{\partial \tilde{x}^{i}}= \nabla_{l}(d\phi(l)),$$
where $\phi^{i}_{ab}= \frac{\partial^{2}\phi^{i}}{\partial x^{a}\partial x^{b}}$, $\bar{G}^{a}$ and $G^{i}$ are the geodesic coefficients of $(M,\bar{F})$ and $(N,F)$ respectively and $h$ is the \emph{normal curvature vector field}, which is in $\mathcal{N}M$ for any given $y$.

The other one is to induce the metric by a unit normal vector field $\xi\in\mathcal{N}^{0}M$, we can define a Riemannian metric $\hat{g}:=\phi^{*}g_{\xi}$. Note that $\xi$ may be locally defined, so the following results hold in a neighborhood of each point on $M$. We call $(M,\hat{g})$ an
\textit{anisotropic submanifold} in $(N,F)$ to distinguish it from an isometric immersion submanifold $(M,\overline{F})$. For any $X\in T_xM$, the \emph{shape operator}~${A}_{\mathbf{\xi}}:T_xM\rightarrow T_xM$ is defined by
\begin{equation*}\label{0.1}
{A}_{\xi}X=-\left(\nabla^{\xi}_{X}\xi\right)^{\top}_{g_{\xi}}.
\end{equation*}
Since
\begin{align}\label{00000}
g_{\xi}(A_{\xi}X,Y)=g_{\xi}(X,A_{\xi}Y),\ \ \ \ \forall X,Y\in T_xM,
\end{align}
we define the eigenvalues of ${A}_{\xi}$. We call $k_{1}, k_{2}, \cdots, k_{m}$ as the principal curvatures of $M$ with respect to $\xi$. If $k_{1} = k_{2} = \cdots = k_{m}$ for any $\xi$, we call $M$ an \textit{anisotropic totally umbilic submanifold} in $(N,F)$. Denote $\hat{H}_{\xi}=\frac{1}{m}\sum\limits_{a=1}^{m}k_{a}$ as the \textit{anisotropic mean curvature} of $M$. If $\hat{H}_{\xi}=0$ for every $\xi$, we call $M$ \textit{anisotropic minimal}. For $\xi(x)\in\mathcal{N}^{0}_{x}M$, let $\xi_{1}$ and $\xi_{2}$ be two different normal vector fields extending by $\xi(x)$. Set $X,Y\in T_{x}M$. We have
\begin{align}\label{0000}
g_{\xi_{1}}(A_{\xi_{1}},Y)|_{x}=-g_{\xi_{1}}(\nabla^{\xi_{1}}_{X}\xi_{1},\tilde{Y})|_{x}=g_{\xi_{1}}(\xi_{1},\nabla^{\xi_{1}}_{X}\tilde{Y})|_{x}=g_{\xi_{2}}(\xi_{2},\nabla^{\xi_{2}}_{X}\tilde{Y})|_{x}=g_{\xi_{2}}(A_{\xi_{2}},Y)|_{x},
\end{align}
where $\tilde{Y}$ is an extension of $Y$. Hence, $A_{\xi}$ is only related to $\xi(x)$ and is not related to its extending vector field.

Let $\rho_{r}$ be the elementary symmetric function of the principal curvatures $k_{1},\cdots,k_{m}$, i.e., $\rho_{r}:=\sum\limits_{i_{1}<\dots<i_{r}}k_{i_{1}}\cdots k_{i_{r}}$ for $1\leq r\leq m$. Set $\rho_{0}=0$. Then the \textit{$r$-th mean curvature} $\hat{H}_{\xi,r}$ is defined by $\hat{H}_{\xi, r}=\frac{1}{C^{r}_{m}}\rho_{r}$, where $C^{r}_{m}=\frac{m!}{r!(m-r)!}$. In particular, $\hat{H}_{\xi,1}=\hat{H}_{\xi}$. We simply write $\hat{H}_{r}$ as $\hat{H}_{\xi,r}$.

Particularly, if $m=n-1$, there exists two global unit normal vector fields $\xi_{\pm}$.
If $F$ is reversible, then $\xi_{-}=-\xi$. For a given $\xi$, $\hat{g}=\phi^{*}g_{\xi}$ can be defined globally. Then $(M,\hat{g})$ is called an oriented $\textit{anisotropic hypersurface}$. From \cite{HYS}, $\hat{\nabla}^{\perp}_{X}\xi=0$, where $\hat{\nabla}$ is the induced normal connection on $(M,\hat{g})$. Hence, in this case, we have
$${A}_{\mathbf{\xi}}X=-\nabla^{\mathbf{\xi}}_{X}\mathbf{\xi}.$$
\subsection{Variation of the induced volume}
Let $\phi: M\to (N,F)$ be an immersion and $d\mu_{M}=\sigma(x)dx$ be a volume form induced by arbitrary ways on $M$. Suppose that $\phi_{t}:M\rightarrow (N,F)$ is a smooth variation of $\phi$, which induces a variation vector field $X:=\frac{\partial \phi_{t}}{\partial t}\big|_{t=0}=X^{i}\frac{\partial}{\partial \tilde x^{i}}\in TN$, where $t\in(-\varepsilon,\varepsilon)$. Denote $V_{t}(M):=\int_{M}d\mu_{M_{t}}$, then
\begin{align}\label{3.9}
V'(0)=-\int_{M}\mathcal{H}_{\sigma}(X)d\mu_{M},
\end{align}
where $\mathcal{H}_{\sigma}\in \mathcal{V}(M)$ is called the \emph{$\sigma$-mean curvature form}.

Particularly, if $M$ is a hypersurface in $(N,F)$, We define the \emph{$\sigma$-mean curvature} $H_{\sigma}$ with respect to $d\mu_{M}=\sigma(x)dx$ by
\begin{align}\label{3.9111}
H_{\sigma}:=\mathcal{H}_{\sigma}(\xi).
\end{align}
If $H_{\sigma}=0$, $M$ is called the \emph{$\sigma$-minimal} hypersurface.

There are many ways to induce the volume form of $M$. For example, we can define the volume form by Busemann-Hausdorff volume form or Holmes-Thmopson voloum form on $M$ with respect to $\bar{F}=\phi^{*}{F}$. For the volume form $d\mu_{M}=\sigma_{\bar{F}}dx$ induced by $\bar{F}$, we express the mean curvature $H_{\sigma_{\bar{F}}}$ by $H_{\bar{F}}$ for simplicity.

From \cite{HS}, we immediately have
\begin{prop}
The mean curvature of $(M,\bar{F})$ with respect to Holmes-Thmopson volume form is
\begin{align}
H_{\bar{F}}=\frac{1}{\int_{S_{x}M}\Omega d\bar{\tau}}\int_{S_{x}M}g(h,\xi)\Omega d\bar{\tau},
\end{align}
where $\Omega=\det(\frac{\overline{g}_{ab}}{\bar{F}})$ and $d\bar{\tau}=\sum\limits_{a}(-1)^{a-1}y^{a}dy^{1}\wedge\cdots\wedge \hat{dy^{a}}\wedge\cdots\wedge dy^{m}.$
Then $(M,\bar{F})$ has constant mean curvature $c_{0}$ in $(N,F)$ if and only if the following Euler-Lagrange equation holds,
\begin{align}
\int_{S_{x}M}(g(h,\xi)-c_{0})\Omega d\bar{\tau}=0.
\end{align}
\end{prop}

There exists another way to induce volume form $d\mu_{M}$ by the global unit normal vector $\xi$. Suppose that $d\mu_{N}=\tilde{\sigma}(\tilde x)d\tilde x$
is a volume form on $(N, F)$ (including Busemann-Hausdorff volume form or Holmes-Thmopson volume form). The induced volume form on $M$ can be defined by
\begin{align}
d\mu_{\xi}&=\sigma_{\xi}(x)dx^{1}\wedge\cdots \wedge dx^{n-1}
\nonumber\\
&=\tilde{\sigma}(\phi(x))\phi^{*}(i_{\xi}(d\tilde x^{1}\wedge\cdots \wedge d\tilde x^{n}))~~x\in M,\label{3.26}
\end{align}
where $i_{\xi}$ denotes the inner multiplication with respect to $\xi$. In this case, we express the mean curvature $H_{\sigma_{\xi}}$ by $H_{\xi}$ for simplicity.
Combining Lemma 14.1.2 in \cite{S} and the proof of Lemma 4.3 in \cite{HYS}, we have
\begin{lemm} \label{lem0_11} The anisotropic mean curvature $\hat{H}_{\xi}$ and the mean curvature $H_{\xi}$ satisfy
\begin{align}\label{3.922}
(n-1)\hat{H}_{\xi}=H_{\xi}-S(\xi).
\end{align}
\end{lemm}

\begin{rema}
It is well known that isoparametric hypersurfaces in Finsler manifolds are hypersurfaces with constant mean curvature. Conversely, the result doesn't come true in general. Similar to Riemannian geometry, we will prove that hypersurface with non-zero constant mean curvature is the critical point of the volume preserving variation.
\end{rema}
\section{Volume preserving variation}
\subsection{The definition in Finsler geometry}
Let $(N^{n},F)$ be a Minkowski space and $d\mu_{N}=\tilde{\sigma}(\tilde x)d\tilde x$ be a volume form of $(N,F)$. For Busemann-Hausdorff volume form or Holmes-Thmopson voloum form, $\tilde{\sigma}(\tilde x)$ is a positive constant. So for simplicity, set $d\mu_{N}=d\tilde x$.

Let $\phi_{0}:M\rightarrow N$ be a closed hypersurface with a given unit normal vector $\xi$ and $\phi_{0}(M)=\phi (M)=\partial U$, where $U$ is a compact domain determined by $\phi_{0}(M)$. Considering a smooth variation $\phi_{t}$ such that $\phi_{t}(M)=\partial U_{t}$, where $U_{t}$ is a compact domain determined by $\phi_{t}(M)$. Let $X=\frac{\partial \phi_{t}}{\partial t}\mid_{t=0}$ be the variation vector field, which can be divided into two parts $X=X^{\top}+X^{\bot}$ with respect to $g_{\xi}$. Obviously, $X^{\bot}=\langle X,\xi\rangle\xi$. We have the following lemma,
\begin{lemm} \label{lem0_1}
$\frac{d}{dt}\text{Vol}(U_{t})|_{t=0}=0$ if and only if $\int_{M}\langle X,\xi\rangle d\mu_{\xi}=0$.
\end{lemm}
\proof
Let $\{\partial_{i}\}$ be the natural basis of $M$ and $\phi_{t}=\phi_{t}(x)$. Using divergence theorem, we have
\begin{align}\label{6.11}
\text{Vol}(U_{t})=\frac{1}{n}\int_{U_{t}}\text{div}\tilde x d\mu_{N}=\frac{1}{n}\int_{M}\langle \phi_{t},\xi_{t}\rangle_{g_{\xi_{t}}} d\mu_{\xi_{t}},
\end{align}
and
\begin{align}
\langle \phi_{t},\xi_{t}\rangle_{g_{\xi_{t}}} d\mu_{\xi_{t}}
&=\langle \phi_{t},\xi_{t}\rangle_{g_{\xi_{t}}} i_{\xi_{t}}(\phi_{t})^{*}d\mu_{N}\nonumber\\
&=i_{\phi_{t}^{\bot}}d\mu_{N}((\phi_{t})_{*}\partial_{1}\wedge\cdots\wedge(\phi_{t})_{*}\partial_{n-1})\nonumber\\
&=d\mu_{N}(\phi_{t}\wedge(\phi_{t})_{*}\partial_{1}\wedge\cdots\wedge(\phi_{t})_{*}\partial_{n-1}).\label{6.12}
\end{align}
Set $V(t)=\text{Vol}(U_{t})$, then we have
\begin{align}
V'(0)
&=\frac{1}{n}\int_{M}\frac{\partial}{\partial t}d\mu_{N}(\phi_{t}\wedge(\phi_{t})_{*}\partial_{1}\wedge\cdots\wedge(\phi_{t})_{*}\partial_{n-1})|_{t=0}\nonumber\\
&=\frac{1}{n}\int_{M}d\mu_{N}(X\wedge\phi_{*}\partial_{1}\wedge\cdots\wedge\phi_{*}\partial_{n-1})\nonumber\\
&+\frac{1}{n}\int_{M}\sum_{i=1}^{n-1}d\mu_{N}(\phi\wedge\phi_{*}\partial_{1}\wedge\cdots\wedge \frac{\partial X}{\partial x_{i}}\wedge\cdots\wedge\phi_{*}\partial_{n-1}).\label{6.13}
\end{align}
By a direct calculation,
\begin{align}
\frac{1}{n}&\int_{M}d\mu_{N}(X\wedge\phi_{*}\partial_{1}\wedge\cdots\wedge\phi_{*}\partial_{n-1})\nonumber\\
&=\frac{1}{n}\int_{M}d\mu_{N}(X^\bot\wedge\phi_{*}\partial_{1}\wedge\cdots\wedge\phi_{*}\partial_{n-1})\nonumber\\
&=\frac{1}{n}\int_{M}\langle X,\xi\rangle d\mu_{\xi},\label{6.14}
\end{align}
and
\begin{align}
\frac{1}{n}&\int_{M}\sum_{i=1}^{n-1}d\mu_{N}(\phi\wedge\phi_{*}\partial_{1}\wedge\cdots\wedge \frac{\partial X}{\partial x^{i}}\wedge\cdots\wedge\phi_{*}\partial_{n-1})\nonumber\\
&=\frac{1}{n}\int_{M}\sum_{i=1}^{n-1}\partial_{i}d\mu_{N}(\phi\wedge\phi_{*}\partial_{1}\wedge\cdots\wedge X\wedge\cdots\wedge\phi_{*}\partial_{n-1})\nonumber\\
&-\frac{1}{n}\int_{M}\sum_{i=1}^{n-1}d\mu_{N}(\phi_{*}\partial_{i}\wedge\phi_{*}\partial_{1}\wedge\cdots\wedge X\wedge\cdots\wedge\phi_{*}\partial_{n-1})\nonumber\\
&=\frac{n-1}{n}\int_{M}\langle X,\xi\rangle d\mu_{\xi}.\label{6.15}
\end{align}
From (\ref{6.13}), (\ref{6.14}) and (\ref{6.15}), we have
$$V'(0)=\int_{M}\langle X,\xi\rangle d\mu_{\xi}.$$
\endproof
According to Lemma \ref{lem0_1}, for any closed hypersurface in a Finsler manifold, we give the following definition,
\begin{defi}
Let $M$ be a closed hypersurface in a Finsler manifold $(N,F)$ with any volume form $d\mu_{M}$. If the variation vector field $X$ satisfies
$$\int_{M}\langle X,\xi \rangle d\mu_{M}=0,$$
where $\xi$ is the unit normal vector, then we say that $X$ determines a volume-preserving variation in Finsler manifolds.
\end{defi}
\subsection{Proof of Theorem \ref{thm1_1}}
Set $\eta=\langle X,\xi\rangle$. From (\ref{3.9}),
\begin{align}
V'(0)=-\int_{M}\mathcal{H}(X^{\bot})dV_{M}=-\int_{M}H_{\xi}\eta dV_{M}.\label{3.92}
\end{align}
Set
$$a=\frac{\int_{M}H_{\xi}dV_{M}}{V(0)},$$
where $V(0)=\text{Vol}(M)$. Namely,
\begin{align}
\int_{M}(H_{\xi}-a)dV_{M}=0.\label{3.911}
\end{align}
On the one hand, if $V'(0)=0$, then $X=\eta\xi$ is the volume preserving vector field. Set $\eta=H_{\xi}-a$, from (\ref{3.92}) and (\ref{3.911}), we have
\begin{align}
0=\int_{M}(H_{\xi}-a)\eta dV_{M}=\int_{M}(H_{\xi}-a)^{2}dV_{M}.\label{3.93}
\end{align}
Obviously,
$$H_{\xi}=a.$$
On the other hand, if $H=constant=c(\neq 0)$,
then
\begin{align*}
V'(0)=-\int_{M}H_{\xi}\eta dV_{M}=-c\int_{M}\eta dV_{M}=0.
\end{align*}
That is, $\phi$ is a critical point of any volume preserving variations.
Thus we get the conclusion.
\section{Geometric correspondence of submanifolds by homothetic navigation}
Let $M$ be a submanifold in a Finsler manifold $(N,F)$ and $\tilde{F}$ be the Finsler metric on $N$ defined by the navigation datum $(F,W)$, where $W$ is a homothetic vector field with dilation $c$. $W$ generates a family(a one-parameter local subgroup) of local diffeomorphisms $\psi_{t}$. $W$ is called a homothetic field if it satisfies
$$(\psi_{t}^{\ast}F)(x,y)=F(\psi_{t}(x),\psi_{t}(y))=e^{-2ct}F(x,y),$$
where $x\in N$, $y\in T_{x}N$ and $t\in\mathbb{R}$. The idea of navigation can be described as follows. We suppose that $F$ is a Finsler metric and $W$ is a vector field with $F(x,-W)<1$, then we can define a new Finsler metric $\tilde{F}$ by
$$F(x, y-\tilde{F}(x,y)W)=\tilde{F}(x,y), \ \ \ \ \ \ \forall x\in N,\ \ y\in T_{x}N.$$
We try to find out the relationship of principal curvatures by homothetic navigation.
\begin{lemm}~\cite{HM12}
Let $F=F(x,y)$ be a Finsler metric on a manifold $N$ and let $W$ be a vector field on $N$ with $F(x, -W_{x})<1$. Suppose that $W$ is homothetic with dilation $c$. Let $\tilde{F}=\tilde{F}(x,y)$ denote the Finsler metric on $N$. Then the geodesics of $\tilde{F}$ are given by $\psi_{t}(\gamma(a(t)))$, where $\psi_{t}$ is the flow of $W$, $\gamma(t)$ is a geodesic of $F$ and $a(t)$ is defined by
\begin{align}
a(t):=\left\{\begin{array}{ll}
\frac{e^{2ct}-1}{2c},\ \ \ \ if\ \ c\neq0,\\
t,\ \ \ \ \ \ \ \ \ \ if\ \  c=0.
\end{array}\right.
\end{align}
\end{lemm}
\begin{lemm}\cite{XM1}
For any orthogonal Jacobi field $J(t)$ along the unit speed geodesic $\gamma(t)$ for the metric $F$,
\begin{align}\label{4.66}
\tilde{J}(t)=(\psi_{t})_{*}(J(\frac{e^{2ct}-1}{2c}))
\end{align}
is an orthogonal Jacobi field along the unit speed geodesic $\tilde{\gamma}(t)$ for the metric $\tilde{F}$. Conversely, any orthogonal Jacobi field $\tilde{J}(t)$ along $\tilde{\gamma}(t)$ for the metric $\tilde{F}$ can be presented by (\ref{4.66}) for some orthogonal Jacobi field $J(t)$ along $\gamma(t)$ for the metric $F$.
\end{lemm}
\begin{lemm}\cite{HD}
Let $\nu\in \mathcal{V}^{0}_{F^{*}}(M)$ and $\xi=\mathcal{L}^{-1}_{F}(\nu)\in\mathcal{N}^{0}_{F}(M)$ be unit normal vectors with respect to $F^{*}$ and $F$, respectively. Then $\frac{\nu}{\tilde{F}^{*}(\nu)}\in \mathcal{V}^{0}_{\tilde{F}^{*}}(M)$ and $\tilde{\xi}=\mathcal{L}^{-1}_{\tilde{F}}(\frac{\nu}{\tilde{F}^{*}(\nu)})\in\mathcal{N}^{0}_{\tilde{F}}(M)$ satisfies
\begin{align}\label{0}
\tilde{\xi}=\xi+W.
\end{align}
\end{lemm}

Using the lemmas above, we have the following theorem,
\begin{theo} \label{thm3.1}
Let $M$ be an anisotropic submanifold in a Finsler manifold $(N,F)$ and $\tilde{F}$ be the Finsler metric on $N$ defined by homothetic navigation with dilation $c$. Let $\lambda$ and $\tilde{\lambda}$ be the principal curvature of $M$ with respect to the unit normal vector $\xi$ in $(N,F)$ and the principal curvature of $M$ with respect to the unit normal vector $\tilde{\xi}$ in $(N,\tilde{F})$, respectively. Then
\begin{align}\label{4.6611}
\tilde{\lambda}=\lambda+c.
\end{align}
\end{theo}
\proof
Set $\tilde{t}=\frac{e^{2ct}-1}{2c}$ and we have
\begin{align}\label{6.99}
\tilde{\gamma}(t)=\psi_{t}(\gamma(\tilde{t}))~~~~~\tilde{J}(t)=(\psi_{t})_{*}(J(\tilde{t})).
\end{align}

From (\ref{6.99}) and Lemma 5.2 in \cite{XM1}, we have
\begin{align}\label{11.1}
\tilde{g}_{\dot{\tilde{\gamma}}}(\tilde{J}_{a}(t),\tilde{J}_{b}(t))=\frac{1}{c_{0}+1}e^{-2ct}g_{\dot{\gamma}}(J_{a}(\tilde{t}),J_{b}(\tilde{t})),
\end{align}
where $c_{0}$ is some real constant. Differentiating both sides along $\tilde{\gamma}$ with respect to $t$. Then the left side of (\ref{11.1}) becomes
\begin{align}\label{11.2}
\frac{d}{dt}\tilde{g}_{\dot{\tilde{\gamma}}}(\tilde{J}_{a}(t),\tilde{J}_{b}(t))
=\tilde{g}_{\dot{\tilde{\gamma}}}(\nabla_{\dot{\tilde{\gamma}}(t)}^{\dot{\tilde{\gamma}}(t)}\tilde{J}_{a}(t),\tilde{J}_{b}(t))
+\tilde{g}_{\dot{\tilde{\gamma}}}(\tilde{J}_{a}(t),\nabla_{\dot{\tilde{\gamma}}(t)}^{\dot{\tilde{\gamma}}(t)}\tilde{J}_{b}(t)).
\end{align}
Let $\Phi_{a}$ be the smooth variation of the geodesic such that $\dot{\tilde{\gamma}}(t)=\Phi_{a*}\frac{\partial}{\partial t}|_{s=0}$ and $\tilde{J}_{a}(t)=\Phi_{a*}\frac{\partial}{\partial s}|_{s=0}$, then we get
\begin{align}\label{11.3}
0=\Phi_{a*}[\frac{\partial}{\partial t},\frac{\partial}{\partial s}]|_{s=0}=[\Phi_{a*}\frac{\partial}{\partial t},\Phi_{a*}\frac{\partial}{\partial s}]|_{s=0}=\nabla_{\dot{\tilde{\gamma}}}^{\dot{\tilde{\gamma}}}\tilde{J}_{a}
-\nabla_{\tilde{J}_{a}}^{\dot{\tilde{\gamma}}}\Phi_{a*}\frac{\partial}{\partial t}|_{s=0},
\end{align}
Notice that $\Phi_{a*}\frac{\partial}{\partial s}$ is the tangent vector field of $M_{t}$ and $\Phi_{a*}\frac{\partial}{\partial t}$ is the normal vector field of $M_{t}$. From (\ref{0000}), $A_{\dot{\tilde{\gamma}}}$ is not related to the extension of $\dot{\tilde{\gamma}}$. Hence, combine (\ref{00000}), (\ref{11.2}) simplifies to
\begin{align}\label{11.4}
\frac{d}{dt}\tilde{g}_{\dot{\tilde{\gamma}}}(\tilde{J}_{a}(t),\tilde{J}_{b}(t))
&=\tilde{g}_{\dot{\tilde{\gamma}}}(\nabla_{\tilde{J}_{a}}^{\dot{\tilde{\gamma}}}\Phi_{a*}\frac{\partial}{\partial t}|_{s=0},\tilde{J}_{b}(t))
+\tilde{g}_{\dot{\tilde{\gamma}}}(\tilde{J}_{a}(t),\nabla_{\tilde{J}_{b}}^{\dot{\tilde{\gamma}}}\Phi_{a*}\frac{\partial}{\partial t}|_{s=0})\nonumber\\
&=\tilde{g}_{\dot{\tilde{\gamma}}}(-A_{\dot{\tilde{\gamma}}}\tilde{J}_{a}(t),\tilde{J}_{b}(t))
+\tilde{g}_{\dot{\tilde{\gamma}}}(\tilde{J}_{a}(t),-A_{\dot{\tilde{\gamma}}}\tilde{J}_{b}(t))\nonumber\\
&=-2\tilde{g}_{\dot{\tilde{\gamma}}}(A_{\dot{\tilde{\gamma}}}\tilde{J}_{a}(t),\tilde{J}_{b}(t)).
\end{align}
where $\tilde{A}_{\dot{\tilde{\gamma}}}$ is the shape operator of $M_{t}$ with respect to $\dot{\tilde{\gamma}}$. Similarly, the right side of (\ref{11.1}) becomes
\begin{align}\label{11.6}
\frac{d}{dt}(e^{-2ct}g_{\dot{\gamma}}(J_{a}(\tilde{t}),J_{b}(\tilde{t})))
=-2cg_{\dot{\gamma}}(J_{a}(\tilde{t}),J_{b}(\tilde{t})))-2g_{\dot{\gamma}}(A_{\dot{\gamma}}J_{a}(\tilde{t}),J_{b}(\tilde{t})).
\end{align}
From (\ref{11.4}) and (\ref{11.6}), we have
\begin{align}\label{111111}
\tilde{g}_{\dot{\tilde{\gamma}}}(A_{\dot{\tilde{\gamma}}}\tilde{J}_{a}(t),\tilde{J}_{b}(t))=\frac{1}{c_{0}+1}(cg_{\dot{\gamma}}(J_{a}(\tilde{t}),J_{b}(\tilde{t}))+g_{\dot{\gamma}}(A_{\dot{\gamma}}J_{a}(\tilde{t}),J_{b}(\tilde{t})))
\end{align}
Let $\{e_{a}\}=\{J_{a}(0)\}$ be the unit orthogonal basis with respect to $g_{\xi}$ and $A_{\xi}e_{a}=\lambda_{a}e_{a}$, where $A_{\xi}=A_{\gamma(0)}$. Set $\tilde{A}_{\tilde{\xi}}e_{a}=h^{d}_{a}e_{d}$, where $\tilde{A}_{\tilde{\xi}}=A_{\dot{\tilde{\gamma}}(0)}$. From (\ref{111111}),
$$\tilde{g}_{\tilde{\xi}}(\tilde{A}_{\tilde{\xi}}e_{a},e_{b})=\tilde{g}_{\tilde{\xi}}(h^{d}_{a}e_{d},e_{b})=h^{d}_{a}\tilde{g}_{\tilde{\xi}}(e_{d},e_{b})=\frac{1}{c_{0}+1}(c+\lambda_{a})\delta_{ab}.$$
From Lemma 3.1 and Lemma 3.3 in \cite{XM1},
$$\tilde{g}_{\tilde{\xi}}(e_{d},e_{b})=\frac{1}{c_{0}+1}\delta_{bd}.$$
Hence,
\begin{align}\label{11.10}
h^{b}_{a}=(c+\lambda_{a})\delta_{ab}.
\end{align}
From (\ref{11.10}), $\tilde{A}_{\tilde{\xi}}e_{a}=(c+\lambda_{a})e_{a}$. Hence, $e_{a}$ is still the principal direction with respect to $\tilde{g}_{\tilde{\xi}}$ and the corresponding principal curvature $\tilde{\lambda}_{a}=c+\lambda_{a}$. This complete the proof of Theorem \ref{thm3.1}.
\endproof
\begin{coro}\label{11.821}
Let $M$ be an anisotropic umbilical submanifold in the Finsler manifold $(N,F)$ and $\tilde{F}$ be the Finsler metric on $N$ defined by homothetic navigation. Then $M$ is also an umbilical submanifold in the Finsler manifold $(N,\tilde{F})$.
\end{coro}

From (\ref{4.6611}), we denote $\hat{\tilde{H}}_{\tilde{\xi}}$ and $\hat{H}_{\xi}$ as the anisotropic mean curvature of $M$ in $(N,\tilde{F})$ and $(N,F)$, respectively. Then $\hat{\tilde{H}}_{\tilde{\xi}}$ and $\hat{H}_{\xi}$ satisfy
\begin{align}\label{00}
\hat{\tilde{H}}_{\tilde{\xi}}=c+\hat{H}_{\xi}.
\end{align}
If $c=0$, then $\hat{\tilde{H}}_{r}=\hat{H}_{r}$ for any $r = 1,\ldots, m$.
Then we have,
\begin{coro}\label{11.831}
Let $M$ be an anisotropic minimal submanifold in the Finsler manifold $(N,F)$ and $\tilde{F}$ be the Finsler metric on $N$ defined by homothetic navigation. If the navigation vector field $W$ is a Killing vector field, then $M$ is also a minimal submanifold in the Finsler manifold $(N,\tilde{F})$.
\end{coro}
\begin{coro}\label{11.841}
Let $M$ be a submanifold in the Finsler manifold $(N,F)$ and $\tilde{F}$ be the Finsler metric on $N$ defined by homothetic navigation $(F,W)$. Then the anisotropic mean curvature of $M$ with respect to $F$ is constant if and only if the anisotropic mean curvature of $M$ with respect to $\tilde{F}$ is constant. Especially when $W$ is a Killing vector field, for some $r = 1,\ldots, m$, the $r$-th mean curvature of $M$ with respect to $F$ is constant if and only if that with respect to $\tilde{F}$ is constant.
\end{coro}

From Corollary \ref{11.841}, we get Theorem \ref{thm1_2} immediately. From Theorem 1.3 in \cite{XM1}, Lemma \ref{lem0_11} and (\ref{00}), we have
\begin{coro}
Let $\tilde{H}_{\tilde{\xi}}$ and $H_{\xi}$ be the mean curvature of $M$ in $(N,\tilde{F})$ and $(N,F)$, respectively. Then $\tilde{H}_{\tilde{\xi}}$ and $H_{\xi}$ satisfy
\begin{align*}
\tilde{H}_{\tilde{\xi}}=2nc+H_{\xi}.
\end{align*}
\end{coro}
In fact, from (\ref{111111}), we obtain
\begin{theo}
Let $M$ be an anisotropic hypersurface in a Finsler manifold $(N,F)$ and $\tilde{F}$ be the Finsler metric on $N$ defined by the navigation datum $(F,W)$ with dilation $c$, where $W$ generates a family of local diffeomorphisms $\psi_{t}$. Let $M_{\tilde{t}}$ and $\tilde{M}_{t}$ be the sufficiently close hypersurface of $M$ with respect to $F$ and $\tilde{F}$, respectively. If $\dot{\gamma}(\tilde{t})\in\mathcal{N}^{0}M_{\tilde{t}}$, then for $x\in M_{\tilde{t}}$ and $\tilde{x}=\psi_{t}x\in\tilde{M}_{t}$, $\dot{\tilde{\gamma}}(t)|_{\tilde{x}}=\psi_{t*}\dot{\gamma}(\tilde{t})|_{x}+W(\tilde{x})\in\mathcal{N}^{0}_{\tilde{x}}\tilde{M}_{t}$. Let $\lambda(\tilde{t})$ and $\tilde{\lambda}(t)$ be the principal curvature of $x$ and $\tilde{x}$ with respect to $\dot{\gamma}(\tilde{t})$ and $\dot{\tilde{\gamma}}(t)$, respectively. Then
\begin{align}\label{1}
\tilde{\lambda}(t)=e^{2ct}(\lambda(\tilde{t})+c).
\end{align}
\end{theo}
\proof
Let $\{J_{a}(\tilde{t}_{0})\}$ be the unit orthogonal principal vector with respect to $g_{\dot{\gamma}}$ and $\lambda_{a}(\tilde{t}_{0})$ be the corresponding principal curvature with respect to $\dot{\gamma}$, where $t_{0}\in(-\varepsilon,\varepsilon)$. Namely, $A_{\dot{\gamma}}J_{a}(\tilde{t}_{0})=\lambda_{a}(\tilde{t}_{0})J_{a}(\tilde{t}_{0})$. Take $\tilde{A}_{\dot{\tilde{\gamma}}}\tilde{J}_{a}=h^{d}_{a}\tilde{J}_{d}$, where $\tilde{J}_{a}(t)=\psi_{t*}(J_{a}(\tilde{t}))$. From (\ref{11.1}), the left side of (\ref{111111}) becomes
\begin{align}\label{2}
g_{\dot{\tilde{\gamma}}}(\tilde{A}_{\dot{\tilde{\gamma}}}\tilde{J}_{a}(t),\tilde{J}_{b}(t))=g_{\dot{\tilde{\gamma}}}(h^{d}_{a}\tilde{J}_{d}(t),\tilde{J}_{b}(t))=h^{d}_{a}\frac{1}{c_{0}+1}e^{-2ct}g_{\dot{\gamma}}(J_{d}(\tilde{t}),J_{b}(\tilde{t})).
\end{align}
Hence,
\begin{align}\label{3}
g_{\dot{\tilde{\gamma}}(t_{0})}(\tilde{A}_{\dot{\tilde{\gamma}}(t_{0})}\tilde{J}_{a}(t_{0}),\tilde{J}_{b}(t_{0}))=h^{b}_{a}(t_{0})\frac{1}{c_{0}+1}e^{-2ct_{0}}.
\end{align}
Consider the right side of (\ref{111111}),
\begin{align}\label{4}
cg_{\dot{\gamma}(\tilde{t}_{0})}(J_{a}(\tilde{t}_{0}),J_{b}(\tilde{t}_{0}))+g_{\dot{\gamma}(\tilde{t}_{0})}(A_{\dot{\gamma}(\tilde{t}_{0})}J_{a}(\tilde{t}_{0}),J_{b}(\tilde{t}_{0}))=(c+\lambda_{a}(\tilde{t}_{0}))\delta_{ab}.
\end{align}
Combine (\ref{3}) and (\ref{4}),
$$h^{b}_{a}(t_{0})=e^{2ct_{0}}(c+\lambda_{a}(\tilde{t}_{0}))\delta_{ab}.$$
Hence, we have
$$\tilde{A}_{\dot{\tilde{\gamma}}(t_{0})}\tilde{J}_{a}(t_{0})=e^{2ct_{0}}(c+\lambda_{a}(\tilde{t}_{0}))\tilde{J}_{a}(t_{0}).$$
That is, for fixed $x=\gamma(\tilde{t})\in M_{\tilde{t}}$ with respect to $F$ and $\tilde{x}=\tilde{r}(t)\in\tilde{M}_{t}=\psi_{t}M_{\tilde{t}}$ with respect to $\tilde{F}$, if $J_{a}(\tilde{t})$ is the principal vector with respect to $g_{\dot{\gamma}}$, then $\tilde{J}_{a}(t)$ is the principal vector with respect to $g_{\dot{\tilde{\gamma}}}$ and $\tilde{\lambda}_{a}(t)=e^{2ct}(c+\lambda_{a}(\tilde{t}))$.
\endproof

\begin{coro}
Let $M$ be an anisotropic hypersurface in the Finsler manifold $(N,F)$ and $\tilde{F}$ be the Finsler metric on $N$ defined by homothetic navigation datum $(F,W)$, where $W$ generates a family of local diffeomorphisms $\psi_{t}$. Let $M_{\tilde{t}}$ and $\tilde{M}_{t}=\psi_{t}M_{\tilde{t}}$ be the sufficiently close hypersurface of $M$ with respect to $F$ and $\tilde{F}$, respectively. Then $M_{\tilde{t}}$ is an anisotropic umbilical submanifold in $(N,F)$ if and only if $\tilde{M}_{t}$ is an anisotropic umbilical submanifold in $(N,\tilde{F})$.
\end{coro}

From (\ref{1}), we denote $\hat{\tilde{H}}_{\dot{\tilde{\gamma}}}$ and $\hat{H}_{\dot{\gamma}}$ as the mean curvature of $M_{\tilde{t}}$ and $\tilde{M}_{t}$ in $(N,\tilde{F})$ and $(N,F)$, respectively. Then $\hat{\tilde{H}}_{\dot{\tilde{\gamma}}}$ and $\hat{H}_{\dot{\gamma}}$ satisfy
\begin{align}\label{X}
\hat{\tilde{H}}_{\dot{\tilde{\gamma}}}(t)|_{\tilde{x}}=e^{2ct}(c+\hat{H}_{\dot{\gamma}}(\tilde{t})|_{x}).
\end{align}
If $c=0$, then $\hat{\tilde{H}}_{r}=\hat{H}_{r}$ for any $r = 1,\ldots, m$.
Then we have,
\begin{coro}
Let $M$ be an anisotropic hypersurface in the Finsler manifold $(N,F)$ and $\tilde{F}$ be the Finsler metric on $N$ defined by homothetic navigation datum $(F,W)$, where $W$ generates a family of local diffeomorphisms $\psi_{t}$. Let $M_{\tilde{t}}$ and $\tilde{M}_{t}=\psi_{t}M_{\tilde{t}}$ be the sufficiently close hypersurface of $M$ with respect to $F$ and $\tilde{F}$, respectively. If $W$ is a Killing vector field, then $M_{\tilde{t}}$ is an anisotropic minimal submanifold in $(N,F)$ if and only if $\tilde{M}_{t}$ is an anisotropic minimal submanifold in $(N,\tilde{F})$.
\end{coro}
\begin{coro}\label{XY}
Let $M$ be an anisotropic hypersurface in the Finsler manifold $(N,F)$ and $\tilde{F}$ be the Finsler metric on $N$ defined by homothetic navigation datum $(F,W)$, where $W$ generates a family of local diffeomorphisms $\psi_{t}$. Let $M_{\tilde{t}}$ and $\tilde{M}_{t}=\psi_{t}M_{\tilde{t}}$ be the sufficiently close hypersurface of $M$ with respect to $F$ and $\tilde{F}$, respectively. Then the anisotropic mean curvature of $M_{\tilde{t}}$ with respect to $F$ is constant if and only if the anisotropic mean curvature of $\tilde{M}_{t}$ with respect to $\tilde{F}$ is constant. Especially when $W$ is a Killing vector field, for some $r = 1,\ldots, m$, the $r$-th mean curvature of $M_{\tilde{t}}$ with respect to $F$ is constant if and only if that with respect to $\tilde{F}$ is constant.
\end{coro}
\begin{rema}
From Corollary \ref{XY} and Definition 2.2 in \cite{DC}, we can prove $M$ is isoparametric in $(N,F)$ if and only if $M$ is isoparametric in $(N,\tilde{F})$.
\end{rema}
From the proof of Theorem 1.3 in \cite{XM1}, denote $\mathbf{S}^{F}(\gamma(\tilde{t}_{0}),\dot{\gamma}(\tilde{t}_{0}))$ and $\tilde{\mathbf{S}}^{\tilde{F}}(\tilde{\gamma}(t_{0}),\dot{\tilde{\gamma}}(t_{0}))$ as the $\mathbf{S}-$curvature in $(N,F)$ and $(N,\tilde{F})$. Then we have
\begin{align}\label{Z}
\tilde{\mathbf{S}}^{\tilde{F}}(\tilde{\gamma}(t_{0}),\dot{\tilde{\gamma}}(t_{0}))=e^{2ct_{0}}\mathbf{S}^{F}(\gamma(\tilde{t}_{0}),\dot{\gamma}(\tilde{t}_{0}))+c(n+1).
\end{align}
Combine Lemma \ref{lem0_11} and (\ref{X}), denote $\tilde{H}_{\dot{\tilde{\gamma}}}$ and $H_{\dot{\gamma}}$ as the mean curvature of $M_{\tilde{t}}$ and $\tilde{M}_{t}$ in $(N,\tilde{F})$ and $(N,F)$, respectively. Then $\tilde{H}_{\dot{\tilde{\gamma}}}$ and $H_{\dot{\gamma}}$ satisfy
\begin{align}\label{Y}
\tilde{H}_{\dot{\tilde{\gamma}}}=e^{2ct}H_{\dot{\gamma}}+c(n+1+e^{2ct}(n-1)).
\end{align}
\begin{rema}
From (\ref{Y}) and Remark 2.3 in \cite{DC}, we have $M$ is $d\mu$-isoparametric in $(N,F)$ if and only if $M$ is $d\mu$-isoparametric in $(N,\tilde{F})$, which proves Theorem 1.5 in \cite{XM1} in a different way.
\end{rema}
\section{Heintze-Karcher Type Inequality}
As we all know, the Minkowski space in Finsler geometry is the natural generalization of the Euclidean space. In this section, we will find out the geometric correspondence of hypersurfaces in Euclidean spaces and Minkowski spaces, based on which we give a Heintze-Karcher Type Inequality in Minkowski spaces.
\subsection{Anisotropic mean curvature in Euclidean spaces}
Recall the definition of anisotropic mean curvatures of hypersurfaces in Euclidean spaces. Consider that $M$ is a closed orientable hypersurface, which is  embedded in $\mathbb{R}^{n}$. Let $F^{*}:\mathbb{S}^{n-1}\rightarrow \mathbb{R}^{+}$ be a smooth positive function satisfying the following convexity condition
$$A_{F^{*}}:=(\nabla^{\mathbb{S}}\nabla^{\mathbb{S}}F^{*}+F^{*}Id_{T_{x}\mathbb{S}^{n-1}})>0$$
for any $x\in\mathbb{S}^{n-1}$. $F^{*}$ can be extended to a 1-homogenous funtion on $\mathbb{R}^{n}$ by
$$F^{*}(x)=F^{*}(\frac{x}{|x|}),\ \ x\in\mathbb{R}^{n}\ \ \textmd{and}\ \ F^{*}(0)=0.$$
In fact, $F^{*}$ is a Minkowski metric, which its dual metric is $F$. Now let $\bar{\nu}$ be the unit normal vector with respect to the Euclidean metric and $\overline{\nu}:M\rightarrow\mathbb{S}^{n-1}$ denote its Gauss map. Then the anisotropic surface energy of $x$ is defined as
$$\mathcal{F}(x)=\int_{M}F^{*}(\overline{\nu})d\mu_{M}.$$
The anisotropic normal vector is defined by
$$\nu_{F}(x):=F^{*}(\bar{\nu}(x))\bar{\nu}(x)+\nabla^{\mathbb{S}}F^{*}(\bar{\nu}(x)).$$
From \cite{MX}, $d\nu_{F}=A_{F^{*}}\circ d\bar{\nu}$. Denote $S_{F}:=-d\nu_{F}$, which is called the \emph{anisotropic Weingarten operator} in a Euclidean space. Define $H_{F}:=tr S_{F}$, which is called the \emph{anisotropic mean curvature} of $M$ in a Euclidean space. The \emph{anisotropic principal curvature} $\kappa^{F}_{1},\cdots,\kappa^{F}_{n-1}$ of $M$ is defined as the eigenvalues of $S_{F}$, then $H_{F}=\sum\limits^{n-1}_{i=1}\kappa_{i}^{F}$.
\begin{lemm}\cite{MX} \label{lem5.1}
Let $x:\Sigma\rightarrow\mathbb{R}^{n}$ be a closed hypersurface embedded into the Euclidean space. If the anisotropic mean curvature $H_{F}$ with respect to the inner normal $\overline{\nu}$ is everywhere positive on $\Sigma$, then we have
$$(n-1)\int_{\Sigma}\frac{F\circ\bar{\nu}}{H_{F}}d\mu_{\Sigma}\geq nV(\Omega),$$
where $V(\Omega)$ is the volume of the compact domain $\Omega$ determined by $\Sigma$. Moreover, the equality holds if and only if $\Sigma$ is anisotropic umbilical.
\end{lemm}

\begin{lemm}\cite{MX}\label{lem5.2}
Let $\Sigma$ be a closed oriented hypersurface embedded in
the Euclidean space $\mathbb{R}^{n}$.  If the $r$-th mean curvature of $\Sigma$ is a non-zero constant for some $r = 1,\ldots, n-1$, then $\Sigma$ is the Wulff shape, up to translations and hometheties.
\end{lemm}
\subsection{Heintze-Karcher Type Inequality in Minkowski spaces}
Let $\xi$ and $\bar{\xi}$ be the normal vectors of $M$ in $N$ with respect to a Minkowski metric and a Euclidean metric, respectively. Let $d\mu_{\xi}$ and $d\mu_{\bar{\xi}}$ be the induced volume form on $M$ determined by $\xi$ and $\bar{\xi}$, respectively.
Then we have the following,
\begin{theo} \label{theo5}
Let $(N,F)$ be a Minkowski space, $U$ be the compact domain on $N$ and $M=\partial U$. If the anisotropic mean curvature $\hat{H}_{\xi}$ with respect to the inner normal vector $\xi$ is positive on $M$ everywhere, then we have
\begin{align}\label{1001}
\int_{M}\frac{1}{\hat{H}_{\xi}}d\mu_{\xi}\geq nV(U),
\end{align}
where $V(U)$ is the volume of $U$. Moreover, the equality holds if and only if $M$ is anisotropic umbilical.
\end{theo}
\proof
Let $H_{\xi}$ be the mean curvature of $M$ with respect to $d\mu_{\xi}$, where $\mathbf{\xi}=\mathcal{L}^{-1}(\nu)$ and $\nu\in \mathcal{V}^{0}(N)$. We have
\begin{align}
d\mu_{\xi}
&=i_{\xi}d\mu\nonumber\\
&=d\mu(\xi\wedge\phi_{*}\partial_{1}\wedge\cdots\wedge\phi_{*}\partial_{n-1})\nonumber\\
&=d\mu(\xi^{\bot}\wedge\phi_{*}\partial_{1}\wedge\cdots\wedge\phi_{*}\partial_{n-1}),\label{4.1}
\end{align}
where $\{\partial_{a}\}=\{\frac{\partial}{\partial x^{a}}\}$ denote the natural basis of $M$.

From $\nu=\frac{\bar{\nu}}{F^{*}(\bar{\nu})}$, we have
\begin{align}\label{10.13}
1=\nu(\xi)=\frac{\bar{\nu}(\xi)}{F^{*}(\bar{\nu})}=\frac{\bar{\nu}(\xi)}{F^{*}(\bar{\xi})}.
\end{align}
Thus we get
$$\xi^{\bot}=\langle\bar{\xi},\xi\rangle\bar{\xi}=\bar{\nu}(\xi)\bar{\xi}=F^{*}(\bar{\xi})\bar{\xi}.$$
Then (\ref{4.1}) simplifies to
\begin{align}
d\mu_{\xi}
&=F^{*}(\bar{\nu})d\mu(\bar{\xi}\wedge\phi_{*}\partial_{1}\wedge\cdots\wedge\phi_{*}\partial_{n-1})=F^{*}(\bar{\xi})d\mu_{\bar{\xi}}.\label{4.2}
\end{align}
Integrating both sides of (\ref{4.2}) and considering the first variation, we have
\begin{align}\label{10.9}
\frac{d}{dt}V(t)|_{t=0}=\frac{d}{dt}\mathcal{F}|_{t=0}.
\end{align}
Let $X$ be the variation vector field, then the left side of the (\ref{10.9}) becomes
\begin{align}\label{10.10}
\frac{d}{dt}V(t)|_{t=0}=-\int_{M}\langle H_{\xi}\xi,X\rangle d\mu_{\xi}=-\int_{M}H_{\xi}\nu(X)d\mu_{\xi}.
\end{align}
The right side of the (\ref{10.9}) becomes
\begin{align}\label{10.11}
\frac{d}{dt}\mathcal{F}|_{t=0}=-\int_{M}\langle H_{F}\bar{\xi},X\rangle d\mu_{\bar{\xi}}=-\int_{M}H_{F}\bar{\nu}(X)d\mu_{\bar{\xi}}.
\end{align}
From (\ref{10.13}), (\ref{4.2}), (\ref{10.10}) and (\ref{10.11}), we obtain
\begin{align}\label{10.12}
\int_{M}(H_{\xi}F^{*}(\bar{\xi})\nu(X)-H_{F}\bar{\nu}(X))d\mu_{\bar{\xi}}=0.
\end{align}
Thus we have
$$H_{\xi}F^{*}(\bar{\xi})\nu-H_{F}F^{*}(\bar{\xi})\nu=0,$$
and
\begin{align}\label{11111}
H_{\xi}=H_{F}.
\end{align}
Combining Remark \ref{rema1}, in Minkowski spaces, we have $H_{\xi}=(n-1)\hat{H}_{\xi}$. Hence, From (\ref{11111}) and the Heintze-Karcher Type Inequality in Lemma \ref{lem5.1}, (\ref{1001}) holds.

In addition, from \cite{XC},
$$\nu_{F}=F_{\tau_{i}}^{*}(\bar{\nu})\frac{\partial}{\partial x^{i}}=F^{*}_{\tau_{i}}(\nu)\frac{\partial}{\partial x^{i}}=\mathcal{L}^{-1}(\nu)=\xi.$$
Since $(N,F)$ is a Minkowski space, set $X\in T_{x}M$,
$$S_{F}X=-d\nu_{F}(X)=-\nabla^{\xi}_{X}\xi=A_{\xi}X.$$
Hence, the principal curvature of anisotropic hypersurface $(M,\hat{g})$ in the Minkowski space $(N,F)$ is the same as $\kappa^{F}_{1},\cdots,\kappa^{F}_{n-1}$. Then $(M,\hat{g})$ is anisotropic totally umbilical if the equality holds in the Minkowski space $(N,F)$. This completes the proof.
\endproof
\subsection{Proof of Theorem \ref{thm1_3} }
Using (\ref{00}), we can get (\ref{1111}). Combining Theorem \ref{thm3.1}, we can obtain that the equality of (\ref{1111}) holds if and only if $M$ is anisotropic umbilical with respect to $\tilde{F}$ either.
\subsection{Proof of Theorem \ref{thm1_4} }
Combining the proof of Theorem \ref{theo5}, Corollary \ref{11.841} and Lemma \ref{lem5.2}, we know that $(M,\hat{g})$ is anisotropic totally umbilical in the Minkowski space $(N,F)$. From \cite{HYS}, we can obtain Theorem \ref{thm1_4}.\\

Yali Chen\\
School of Mathematics and Statistics, Anhui Normal University, Wuhu, Anhui, 241000, China\\
E-mail: chenyl@ahnu.edu.cn\\

Qun He \\
School of Mathematical Sciences, Tongji University, Shanghai, 200092, China\\
E-mail: hequn@tongji.edu.cn\\

Yantong Qian\\
School of Mathematical Sciences, Tongji University, Shanghai, 200092, China\\
E-mail: 1810875@tongji.edu.cn
\end{document}